\documentclass[11pt,a4paper]{article}
\usepackage[OT4,OT1,plmath]{polski}
\usepackage[OT2,T1]{fontenc}

\input cyracc.def
\newcommand\cyrfamily{\fontencoding{OT2}\fontfamily{wncyr}
\selectfont\cyracc}
\DeclareTextFontCommand{\textcyr}{\cyrfamily}

\usepackage[polish,english]{babel}

\usepackage[utf8x]{inputenc}
\usepackage{amsmath}
\usepackage{graphicx}
\usepackage{wrapfig}

\graphicspath{{./Figures/},{../Figures/},{./Pictures/}}

\usepackage{booktabs, cellspace}

\usepackage[colorinlistoftodos]{todonotes}
\usepackage{hyperref}
  \hypersetup{
    pdftitle="Design of engineering systems in Polish mines in the second half of the 20th century", 
    pdfauthor=Aneta Antkowiak and Monika Kaczmarz and Krzysztof Szajowski, 
    colorlinks,
    urlcolor=blue,
    filecolor=magenta,
    citecolor=green, 
    linkbordercolor={1 1 1}, 
    citebordercolor={1 1 1},  
    urlbordercolor={ 1 1 1}  
  } 
\RequirePackage[numbers]{natbib}
	
\newif\ifwww
\wwwtrue 
\def\repo{http://wydawnictwa.ptm.org.pl/index.php/matematyka-stosowana/article/viewArticle}  

\newcommand*{\doi}[1]{\href{http://dx.doi.org/#1}{doi: #1}}
\newcommand*{\MR}[1]{\href{http://www.ams.org/mathscinet-getitem?mr=#1&return=pdf}{MR #1}}
\newcommand*{\ZBL}[1]{\href{http://www.zentralblatt-math.org/zmath/en/advanced/?q=an:#1&format=complete}{Zbl #1}}

\title{Design of engineering systems in Polish mines in the third quarter of the 20th century\footnote{The article is based upon  work created in Wroc{\l}aw University of Technology.}
}
\author{Aneta Antkowiak, Monika Kaczmarz and Krzysztof Szajowski}

\begin{document}
\maketitle

\begin{abstract}
Participation of mathematicians in the implementation of economic projects in Poland, in which mathematics-based methods played an important role, happened sporadically in the past. Usually methods  known from  publications  and verified were adapted to solving related problems. The subject of this paper is the cooperation  between mathematicians and engineers in Wrocław in the second half of the twentieth century established in the form of an analysis of the effectiveness of engineering systems used in mining. The results of this cooperation showed that at the design stage of technical systems it is necessary to take into account factors that could not have been rationally controlled before. The need to explain various aspects of future exploitation was a strong motivation for the development of mathematical modeling methods. These methods also opened research topics in the theory of stochastic processes and graph theory. The social aspects of this cooperation are also interesting.
\end{abstract}

\section{Preliminaries.} In the first half of the twentieth century significant advances in mathematical models of random phenomena led to the emergence of a separate area of mathematics using advanced algebra and mathematical methods to solve engineering problems. Although the history of probability is much longer (see the section \ref{ModeleS}), the use of probabilistic methods in technical sciences has not been easy to trace in the past. Research topics were connected to questions raised by professionals in various fields. Omission of  detailed features of the studied phenomena in the mathematical modeling process often led to similar descriptions of different issues. Universal mathematical models  also became a tool used in  case studies. After the Second World War, the mining industry was one of the developing industries in Poland. Existing mines were modernized and    new ones were opened: coal, lignite, copper. The industrial extraction of raw materials required a rational scheme of both the modernization and design of new mining plants (cf. \cite{Cza11:Meaning}, \cite{BatGlaSaj73:Zarys}).Usually is not  possible to base  the construction of such plants on the known, functioning solutions,  due to the unique nature of the place, tasks and objectives. Expenditure is significant and  there is high risk of improper selection of the structure of the operating system and its components.

Electricity production is a very sensitive part of the economy. It is not strange, therefore, that the Polish government, and in particular the management of financial resources, were interested in the effectiveness of investments in this sector. In the early 1960s, the subject met with the interest of the Investment Bank of Poland, which commissioned an analysis of the efficiency of electricity production at the Turoszów power plant to Włodzimierz W. Bojarski, PhD\footnote{\href{http://www.ippt.pan.pl/staff/WłodzimierzBojarski}{Włodzimierz W.~Bojarski}  \href{http://ww2.senat.pl/k1/senat/Senator/bojarski.htm}{was in the Polish Senate of the first term}.}, who at that time conducted scientific research on a similar subject (reports on this subject in \emph{Życie Warszawy}) in November 1964 (v. \cite{ZW1964:Jeszcze26xi1964}.) By using statistical methods it was found that lignite coal mines in Lower Silesia have many failures. As a result, it was concluded that increasing the reliability of the coal supply system from the mine to the power plant  offers creating
 potential reserves, which may be a way to increase the efficiency of the energy sector. A detailed mathematical model of the engineering system related to transport in the mine could lead  to creating solutions improving delivery of  fuel to the plant. The method of mathematical modeling of the transport system in the opencast mines of Lower Silesia, was developed by Poltegor, which strengthened its potential by engaging the academic staff of Lower Silesia universities with its difficult issues. Available documents of this collaboration are research reports and research articles (v. \cite{Gla64:Process}). New results by mathematicians in Wrocław were presented at the seminar in the library of the Institute of Mathematics of the University of Wroclaw, at that time located in the building of the Faculty of Sanitary Engineering of Wrocław University of Technology\footnote{This building, belonging to Wrocław University of Technology, hosted mathematicians employed at  the University of Wroclaw along with those employed at Wrocław University of Technology.}. Hugo Steinhaus wrote in \cite[Ch. 6, p.~321--322]{HugoSte2016:Vol2}

\begin{quote}
DECEMBER 12, 1963. {\footnotesize ,,\emph{We live in a permanent catastrophe}'', said Karl Kraus about Austria fifty years ago. His words apply today to our past Polish reality: Any day you choose will look something like this: no meat available except expensive tenderloin, none of the better kinds of sausage or salami, no wieners, no lemons, and no sugar. At the post office there are no ordinary or airmail envelopes — and when all of a sudden some turn up they can’t be sealed because the strip of glue is only about a millimeter wide — and no letter paper. Children sit in their classrooms in their winter coats because of the shortage in heating fuel, and it’s the same story at the university.}

{\footnotesize The newspapers strive to distract people’s attention from these shortages with childish exclamations about this, that, or the other largely imaginary achievement. Thus, for instance, they write about the number of megawatts generated by the Turów power plant, yet here I am sitting in the unheated library of the Mathematical Institute of University of Wroclaw listening to a lecture by Gładysz, PhD, a mathematical expert from that same plant, in which he concludes that ,,Turów is doomed''.}
\end{quote}  

Gładysz\footnote{The bio of Stanisław Gładysz (1920--2001) is presented in Appendix~\ref{StGladysz} and was published by Wiadomości Matematyczne~\cite{BycKasKroRom03:Gladysz}.} worked on this issue with Jerzy Battek\footnote{The bio of Jerzy Battek(1927--1991) and the information about his scientific achievement are  presented in Appendix~\ref{JerzyBattek}. See also \href{http://pti.wroc.pl/html/pdf/historiaInformatyki/BattekJerzy_ZHuzar1991.pdf}{text on site  http://pti.wroc.pl}.} and Jan Sajkiewicz\footnote{The biographical sketch  of Jan Sajkiewicz (1928--1995) and the information about his scientific achievement are  presented in Appendix~\ref{JanSajkiewicz}.}. We will devote Chapter \ref{OCT} to this cooperation and its results. The academic implications of this subject, in the form of doctoral dissertations and new fields of study, and the social implications of the participation of mathematicians in solving research problems faced by the industrial sector will be also presented.

Information about research on similar problems in the world at that time reached Poland with a long delay. Famous mathematicians were working on mathematical modeling of random phenomena, and Polish mathematicians were keen on this subject. At this point, it is worthwhile to trace  what results in the theory of stochastic processes and their applications to issues similar to those occurring in the mining industry were most likely to be known in Poland. We refer to the analyzes of works conducted and published by mathematicians working at the Mathematical Institute of the University of Wrocław, such as Bolesław Kopociński, Ilona Kopocińska\footnote{\href{http://www.math.uni.wroc.pl/~ibk/ik/cv.html}{Ilona Kopocińska (1938--2016)}.}, Józef Łukaszewicz\footnote{\href{http://www.ibspan.waw.pl/komisja.statystyki/wspomnienia/Lukaszewicz.pdf}{Józef Łukaszewicz (1927--2013)} v. \cite{Wyl2014:JLukaszewicz}.} and co-workers and in the Mathematical Institute of the Polish Academy of Sciences: Eugeniusz Fidelis\footnote{Eugeniusz Fidelis (1927--2014) v. \cite{Ste2016:Fidelis,Ste2015:Fidelis}.}, Andrzej Wakulicz\footnote{Andrzej Wakulicz (1934--2010) v. \cite{Reg2011:Wakulicz}.}, Ryszard Zieliński\footnote{\href{http://www.ibspan.waw.pl/komisja.statystyki/wspomnienia/Zielinski.pdf}{Ryszard Zieliński (1932--2012) } v. \cite{Zie2013:RZielinski,Nie2012:RZ,Bor2012:RZ,Ryc2012:RZ,WZiel2012:RZ2,WZiel2012:RZ1}. } and co-workers. We will discuss the aspects in chapter \ref{ModeleS}, where we focus on research into Markov processes that are the basis of models used to solve transport problems in mines.

The presented systemic approach to the analysis of technology and structure of engineering systems through the study of the mathematical model allows for conducting important analyzes in isolation from the actual data on the operation of the system, including transport and mining systems before their launch. However, it is not possible to assess the real situation without collecting information about phenomena, events and effects from the existing, operating system. This should be done not only for the assessment of the model, but also for reporting phenomena and managing existing systems. This topic, although very interesting and important, is not addressed in this study. We only mention that in the design of new systems in the 1960s-1980s  such data were usually the only basis for traditional methods of designing engineering systems for new mines. The  methods developed on the basis of mathematical models gained importance at the end of the 20th century, when the verification of project assumptions was given to IT specialists who, using computer simulations, assessed the effects of design assumptions.

\section{A story about mathematics and open pit mines.}
The need to develop the theory of operation of machine systems with continuous transport was noticed in the early sixties of the twentieth century. It was after the first experiments in the work of large opencast lignite coal mines in Poland (,,Turów'', ,,Adamów''), in which belt transport of coal and overburden was applied. Conveyor belts are  pieces of equipment that  have very high reliability. They also have a number of other merits and for this reason  are recommended for wide application. The most sensitive as well as the most expensive item of conveyors is the belt. To keep its durability at high level even in mining operations when tough broken pieces of rock are transported, a~sizer is frequently applied to crush rock before its entrance on a belt. A~problem  considered here is an evaluation of reliability of a system: loading machine - crusher - belt conveying. Today we know that the applied mathematical modeling  by means of Markov processes in this regard is justified in a small number of practical cases (cf. \cite{Cza09:Belt}). This pioneering work was performed by Battek, Gładysz and Sajkiewicz with cooperators from Wrocław University of Technology.  After one year, these models were extended by applying semi-Markov processes (v. \citeauthor{Gra2014:Semi-Markov}(\citeyear{Gra2014:Semi-Markov}), \citeauthor{LisFreDin2010:Multi}(\citeyear{LisFreDin2010:Multi})). Such an approach is greatly valuable because of high level of generalization and  allows one to construct basic reliability measures for the system considered in a general case.

\subsection{Hugo Steinhaus notes on mathematics for industry.}
The problems of the energy industry   were under permanent observation by journalists. In the Summer of 1964 there was a report (v. \cite{ZW1964:Szperkowicz23vii1964}, \ref{AppZW1964:Szperkowicz23vii1964}) by Jerzy Szperkowicz, the special envoy of ,,Życie Warszawy''. He described the efforts of mathematicians working for Poltegor to determine weak points of the engineering system in the Turów power plant and its supplier, the lignite coal mine. Steinhaus \cite[Ch. 6, p.330]{HugoSte2016:Vol2} commented in that article:
\begin{quote}{\footnotesize 
I must mention an event about which our press was very sparing in its coverage. It saved  from collapse the enterprise represented by the Turów II mine and Wroclaw mathematicians. A large team of mine directors and engineers proved unable to deal with the problem of stoppages in the motion of the conveyors transporting lignite coal to big electrical turbines. Even capable, experienced engineers were unable to understand why the halts were so frequent, exceeding their probability estimates. The Investment Bank declared that the breakdowns in the supply of coal decreased the efficiency of the generators below the level at which the production of current ceases to be profitable, and there was even talk of an investigation by the public prosecutor into suspected wrong-doing. But the situation was saved by   Stanisław Gładysz PhD, who, using the theory of stochastic processes, was able to provide an efficient plan for the transportation of coal along the network of conveyors. The regional director of the coal mine stated that had they had Gładysz’s PhD input earlier they would have saved billions in unnecessary investment, and that the overall costs of managing the mine and plant would as a result decrease by 10\%. It follows that Wroclaw contributes more to the mining industry 108 than Warsaw takes out of the treasury for the Workshop of Mathematical Apparatuses of the Polish Academy of Sciences. For his contribution Gładysz, PhD, was made a consultant to the plant with a monthly salary of 2000 złotys--less than the average wage of a Turów miner.}
\end{quote}  

\subsection{,,Życie\,Warszawy''\,about\,mathematics\,for\,industry.}\,A\,few months later the topic return and appeared on the pages of ,,Życie Warszawy'' (v.~\cite{ZW1964:Jeszcze26xi1964}). 
\begin{quote}
,,Życie Warszawy'' November 26, 1964

{\bf More about ,,\emph{mathematics on the outcrop}''\footnote{Translation: Joanna Wasilewska.}.}\\
{\footnotesize In the middle of summer we published a feature by Jerzy Szperkowicz on the use of mathematical methods used to analyze mining machines and conveyors in the Turoszów mine. It concerned the cooperation of science and industry and precisely the application of mathematics in industry. It also discussed the achievements of Wroclaw’s mathematicians, however, the work undertaken earlier on the initiative of the Investment Bank was not mentioned. From there we received a letter with an explanation of the issue. We are printing a part of the letter, apologizing to the sender for delay, which was due to technical reasons.

\emph{Already in 1962, contemporary mathematical probabilistic methods were used to analyze mining machines and conveyor belts in the report by the Investment Bank, which was written by Włodzimierz Bojarski, PhD, Eng.\footnote{See also \href{https://pl.wikipedia.org/wiki/Wlodzimierz_Bojarski}{Włodzimierz Witold Bojarski}.}. The author developed an approximate method for calculating the effective mechanical efficiency, adapted to the approximate nature of the output data, which can be used directly in a design office. Following this method, Bojarski, PhD, did, in the bank, independent complete  calculations of the planned machines and conveyor's systems for the Turów II mine, without the help of mathematical machines.}

\emph{Stanisław Gładysz, PhD, from the Mathematics Department at Wroclaw University of Technology, who was mentioned in the article in ,,Życie Warszawy'',  a reviewer of the work by Włodzimierz Bojarski, PhD, having at his disposal a mathematical machine, specified in a certain way and developed the previously described method. Unfortunately, the specification is currently of a minor value, given the estimated nature of the output data.}

\emph{Both methods, however, allow  a much more rational design of mining machinery systems and conveyors than it has been so far.}

We are gladly publishing this explanation, as it is important that the merits related to this interesting innovation  be clearly highlighted. 

The letter from the Investment Bank also points out that in the case described by J. Szperkowicz, it is  difficult to  talk about substantial savings  as ,,\emph{in relation to the machines, which are currently constructed in open-pit mines, these calculations indicate serious under-investment and the need for costly supplements.}''

Obviously, it does not discredit the connections of science and practice. On the contrary, it shows the need for a close relationship between them.}
\end{quote} 

\subsection{Further discussion in ,,Życie Warszawy''.} This press release, according to the author's suggestion, was  given a riposte. Written  by prof. Gładysz, and supported by detailed explanations of the originality and value of the method proposed by the team which included, among others, Jerzy Battek, Stanisław Gładysz, Jan Sajkiewicz, it was sent to the Dean of the Faculty of Mathematics, Physics and Chemistry of University of Wroclaw. Steinhaus wrote about the creation of this d\'ementi (v. \cite[Ch. 6, p.332]{HugoSte2016:Vol2}):
\begin{quote}{\footnotesize 
I was recently shown a letter from the Investment Bank to ,,Życie Warszawy'' (Warsaw Life) stating that one of the engineers appointed by the bank to look into the situation at Turów II is responsible for the corrections subsequently ascribed to Gładysz PhD. The letter is short and doesn’t quibble in its denial that Gładysz’s report on the problem resulted in a very significant reduction in the amount of capital investment needed to keep the plant running. The letter makes out that although the engineer had the solution, he was pre-empted by Gładysz because the latter had the advantage of access to an electronic computer. (It is slightly difficult to know  what exactly the engineer had prepared, since there is no information on whatever had been done before the arrival of Gładysz on the scene.) However it may be that, this is the first time I have heard of a bank touting itself as a fount of mathematical expertise. I feel I can elaborate on this since a few years ago I was asked to evaluate a report written by the expert sent by the bank to Turów.}
\end{quote}     
Some\,fragments\,of\,this\,reply\,were\,published\,by\,,,Życie\,Warszawy''\,(v.~\cite{ZW1964:Poza11xii1964}, \ref{AppZW1964:Poza11xii1964}), and Steinhaus wrote:
\begin{quote}
DECEMBER 17, 1964\,$\ldots$\\ {\footnotesize
Here’s another story typical of Poland today. Two years ago--as I mentioned above--Turów (a lignite coal mine and power plant) was at an impasse: breakdowns of the conveyor belts reduced the work time of the power plant to 70\%--and it is well known that a large plant ceases to be profitable if breakdowns absorb 30\% of its work time.  This was brought to light by Mr. Bojarski, an expert representing the Investment Bank. Those in charge at the mines turned to Wroclaw University of Technology for help, where the engineering professors declared the mathematical technicalities beyond them, and handed the problem over to the mathematician Stanisław Gładysz, who, together with Battek PhD and others, solved this exceedingly difficult problem. (Several papers on their solution were subsequently published in the periodical \emph{Węgiel Brunatny} (Lignite Coal).) Well, imagine the surprise when not long ago the magazine ,,Życie Warszawy'' printed a communication from the Investment Bank that credit for solving the problem is entirely due to Bojarski, Gładysz’s contribution being essentially insignificant. The Bank’s letter went so far as to hint that Gładysz had plagiarized Bojarski’s work. We all knew at the university that there was not a~jot of truth in this since Bojarski’s paper merely described the factual situation at the mine, and was moreover written in a rather naive way, simply giving the statistics regarding the inefficiency of the system and deducing the loss in production of electricity, whereas Gładysz gave concrete advice on how to organize the network of conveyor belts to ensure their steady motion. Professors Marczewski, Rybarski, Ślebodziński, Urbanik, and I decided to act to put an end to the obfuscation with which the Investment Bank had gone about establishing scientific priorities. We were aided in this by the Bureau for the Project to Mine Lignite Coal, whose director, Schmidt, confirmed the facts of the matter without hesitation, and ,,Życie Warszawy'' has published every argument we supplied in our effort to ensure that truth prevails. Our rectification will in all likelihood reach Minister Szyr, the Chairman of the Committee for Science and Technology.}
\end{quote}

\section{\label{ModeleS}Applications of probability in the second half of the 20th century.}
The topic is considered from the Polish perspective in the 1960s  and takes into account the  methods used by \citeauthor{BatGlaSaj73:Zarys}(\citeyear{BatGlaSaj73:Zarys}). The main tools were restricted to the Markov processes with finite sets of states. Let us recall a few facts from the history of  Markov processes.   

\subsection{Markov properties of random sequences.} In early history of the probability theory (v.~\cite{She74:prehistoryProbab}) the idea of sole dependence of the last result of an experiment appears in various works. Records related to the issue can be found in \cite{Wie2006:TekstyProbab,Wie2013:ProbabStat}. In the text of Jan Śniadecki, ,,The calculation of events and cases of chance''\footnote{Full text is published in \cite{Wie2015:PierwszyTextProbab}.} dated 1790, the author of the outline discusses independent and dependent events (cf. \cite[p.102]{Wie2013:ProbabStat}. Andrey Markov studied \emph{Markov chains} in the early 20th century (v.~\cite[Ch. 13]{She12:HStatistics}, \cite{Gag17:Markov}, \cite{Sen1996:Markov}). At that time the investigation of an extension of independent random variables was motivated, among others, by finding the necessary condition for the weak law of large numbers to hold. Markov~\cite{Mar1906:Extension} motivated by the disagreement with Pavel Nekrasov, who claimed independence was necessary for the weak law of large numbers to hold,  showed that under certain conditions the average outcomes of a~\emph{Markov} chain would converge to a~fixed vector of values, thus proving a~weak law of large numbers without the independence assumption, which had been commonly regarded as a requirement for such mathematical laws to hold (cf. \cite{Hay13:Markov}). Markov later used Markov chains to study the distribution of vowels in ``Eugene Onegin", written by Alexander Pushkin (v. \cite{Mar1913:Onegin}), and proved a central limit theorem for such chains. Sheynin~\cite[Ch. 13.3]{She12:HStatistics} wrote: ,,\emph{4) Markov chains. This term is due to Bernstein~\cite[§16]{Ber1926:dependantes}}. Markov himself (cf. \cite[p.354]{Mar1906:Extension} called them simply chains, but he never wrote about their applications to physical  science (v.~\cite{GriSne97:Introduction}). He issued from a paper by Bruns~\cite{Bru1906:Wahr}, but the prehistory of Markov chains is much richer.''  There are opinions that ,,Questions about dependence $\ldots$ were further raised in a correspondence between $\ldots$  Markov and $\ldots$ Chuprov, which became a starting point of a general theory of stochastic processes'' (cf. \cite[p. 115]{GigSwiPorDasBeaKru89:Empire}, \cite{You1974:Markov}).

In 1912 Poincaré studied Markov chains on finite groups with an aim to study card shuffling (cf.~\cite{Poi1912:shuffling}). Other early uses of Markov chains include a diffusion model introduced by Paul and Tatyana Ehrenfest in \cite{EfrPT1907:diffusion}, and a branching process introduced by Francis Galton and Henry William Watson in 1873\footnote{The history of study \emph{a branching process} started with the problem posed by Francis Galton(1873) in \cite{Gal1873:4001} (v. \cite[p. 386]{Ken1966:Branching}, \cite{Ken1975:Genealogy}, \cite[p. 377]{GriSne97:Introduction}).} , preceding the work of Markov (cf. \cite[p.~377-379]{GriSne97:Introduction} and \cite{Bre1999:Gibbs}). 

After the work of Galton and Watson, it was later revealed that their branching process had been independently discovered and studied around three decades earlier by Irénée-Jules Bienaymé (v. comments and translations in \cite[pp. 117-118]{HeySen1977:Bienayme} and \cite{Sen1996:Markov,Sen1998:Early}). Starting in 1928, Maurice Fréchet became interested in Markov chains, which eventually resulted in him publishing in 1938 a detailed study on Markov chains.
The history of Russian mathematics in the 19th century up to 1930 can be found in \cite{KolYus1978:16century,KolYus1992:16century}. 

Andrei Kolmogorov\footnote{There are many historical papers about A.N. Kolmogorov and his influence on the contemporary state of the probability theory and its application. The early writings on the topic are \cite{Ken1990:Kolmogorov} and \cite{Cra1976:Half}.} developed in the 1930s  a~large part of the early theory of continuous-time Markov processes. He knew Louis Bachelier's 1900s work on fluctuations in the stock market as well as Norbert Wiener's work on Einstein's model of Brownian movement and followed their ideas (v.~\cite{Ken1990:Kolmogorov}). He introduced and studied a~particular set of Markov processes known as diffusion processes, where he derived a~set of differential equations describing the processes (v.~\cite{Sko2005:Basic}). Independent of Kolmogorov's work, Sydney Chapman\footnote{Sydney Chapman was born in Eccles, near Salford in England on 29th of January 1888. He was a British mathematician and geophysicist. His work on the kinetic theory of gases, solar-terrestrial physics, and the Earth's ozone layer inspired a broad range of research over many decades (v.~\cite{Cow1971:Chapman}).  Chapman died in Boulder, Colorado, on 16th of June 1970.} derived in \cite{Cha1928:diffusion} a similar equation. The explanation of the derivation was less mathematically rigorous than Kolmogorov result (v.~\cite{Ber2005:Bachelier}).  The differential equations are now called the Kolmogorov equations (cf.~\cite{And1991:Continuous}) or the Kolmogorov–Chapman equations (v.~\cite{Ken1990:Kolmogorov}). Other mathematicians who contributed significantly to the foundations of Markov processes include William Feller\footnote{Early descriptions of Feller's achievements were given in \cite{Kac1972:Feller}, \cite{Doo1972:Feller}. A recent bio of Feller  is a part of \cite{Fel2015:VolI}, \cite{Fel2015:VolI}. His publication \cite{Fel1950:Introduction} is commonly recognized as one of the great events in mathematics of XX century.}, starting in 1930s, and then later Eugene Dynkin, starting in the 1950s (v.~\cite{Cra1976:Half}).

\subsection{Examples of probabilistic modeling from this period.}
Let us now recall a few spectacular results on the probabilistic applications of that period which have reached Wrocław and stimulated the use of probabilistic models in engineering sciences.

\subsubsection{The storage problems} were subject of research in many research centers. \citeauthor{Gani1957:storage}(\citeyear{Gani1957:storage}) presented a valuable review of work on them. The paper discussed: provisioning, the problem of emptiness and the optimal inventory policy. Important questions concerned  also stationary distributions for finite and infinite dams as well as numerical methods  (v. also \cite{VitThaMur2013:storage} and the monograph  \cite{Mor1959:storage}). 

\subsubsection{Studying the properties of materials} 
\begin{figure}[th!]
\centerline{\includegraphics[width=90mm]{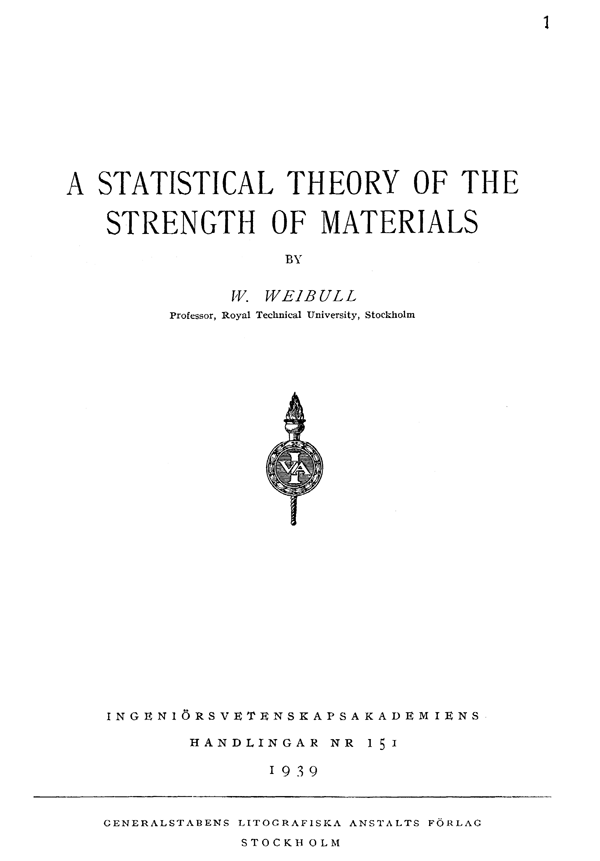}}
\caption{\label{Weibull1939}The first page of the article \cite{Wei1939:Dist}\protect\footnotemark.}
\end{figure}

used in construction works also requires  mathematical models (v. \cite{Irw1957:stress}). Many engineers dealt with this issue. The results of such research  lead to the creation of a mathematical model of statistical material behavior. Waloddi Weibull obtained a unique result in this respect (v. \ref{EHWWeibull}) 
He was interested in modeling the behavior of materials (v. \cite{Wei1927:dynamiska}). His research is of  greatest importance, since it led to the introduction of a probability distribution that describes well various experimental measurements. In 1939 Weibull \cite{Wei1939:Dist} (see Fig.~\ref{Weibull1939}) presented the mathematical model and then he was invited in 1951 to present his most famous result to the American Society of Mechanical Engineers (ASME) on the distribution related to the statistical analysis of many kinds of measurements.The lecture  was followed by a publication with attached discussions \cite{Wei1951:StatDist}. The Weibull distribution is named after him although it was first identified in 1927 by Fréchet \cite{Fre1927:BullPAS} and first applied in 1933 by Rosin \& Rammler\cite{RosRam1933:Laws} (v. \cite{RosRam1927:Laws}) to describe a particle size distribution. Weibull obtained his doctorate from the University of Uppsala in 1932. He was employed in Swedish and German industries as a consulting engineer. 

\footnotetext{\href{http://www.barringer1.com/wa.htm}{Weibull's articles from Dr. Robert B. Abernethy’s Library}.}

Of course, the search for good mathematical models for experimental data continues to this day. In the field of material fatigue modeling, it is worth mentioning the result of Birnbaum${}^{\text{\cite{Woy2001:Birnbaum,Woy1997:Birnbaum}}}$-Saunders \cite{BirSau1969:NewProb,BirSau1969:Stat}, who proposed a distribution also known as the fatigue life distribution, which is used extensively in reliability applications to model failure times.

\subsubsection{Queues} are a phenomenon resulting from the specialization in the scope of services or the way of their implementation, leading in consequence to concentration. When the number of service providers is small and those interested in their service are numerous, it is impossible  that everyone can be satisfied at the same time. There was a need for efficient management of this phenomenon, thus mathematical models of such situations were created. The pioneering work was done by Agner Krarup Erlang (1878--1929) (v. \ref{AKErlang} and \cite{BroHalJen1960:Erlang})

\begin{figure}[th!]
\centerline{\includegraphics[width=9cm]{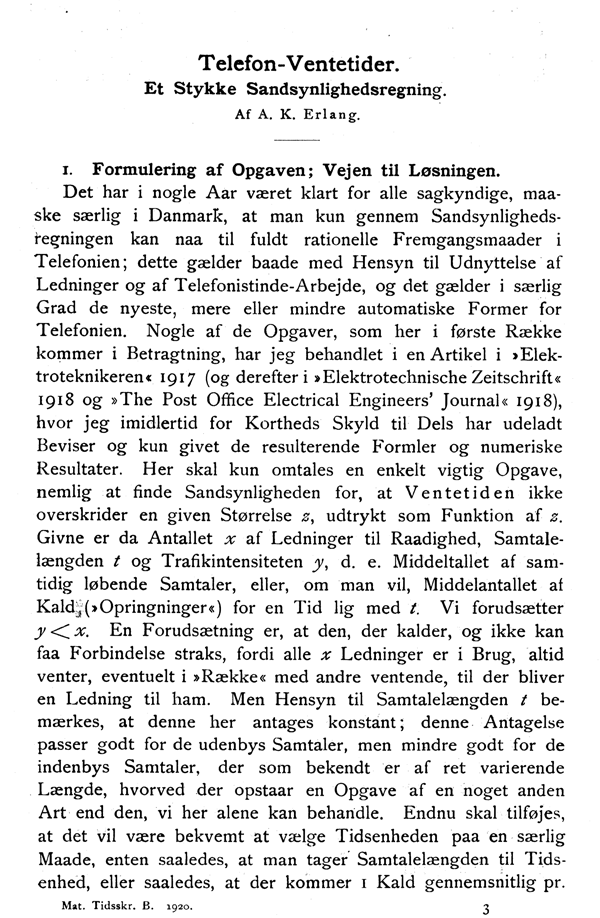}}
\caption{\label{Erlang1920}The first page of the article \cite{Erl1920:Telefon}\protect\footnotemark.}
\end{figure}\footnotetext{\href{http://runeberg.org/}{Project Runeberg}.}

Next, many engineers and mathematicians worked on understanding the phenomena and contributed to the theory of queues: Felix Pollaczek (v. \cite{Coh1981:Pollaczek}), Alexander Iacovlevich Khinchin (v. \cite{Gne1961:Khinchin}), David George Kendall\footnote{The biographical information of D.G. Kendall(1918--2007) can be found in \cite{Kin2009:Kendall}. }. The number of possible queuing variants means that we also have many models. In order to systematize the research, a schematic way of description was introduced. Single queuing nodes are usually described using Kendall's notation in the form A/S/C where A describes the time between arrivals to the queue, S the size of jobs and C the number of servers at the node. Many theorems in queuing theory can be proved by reducing queues to mathematical systems known as Markov chains. Certain elements of modeling queues are analogous to issues related to the description of communication structures in excavations.

\subsubsection{\label{RelModHis}Reliability models.}
An excellent review and sketch of the contemporary problems of the reliability theory and practice can be found in the lecture by Ushakov (v. \cite{Ush2000:Reliability}). It is an overview of general theory. However, such problems are a challenge in the construction of adequate models for the considered engineering constructions or projects. The main problem in such a case is lack of field data. In the case of the engineering equipment for  mining,  additional difficulties are related to uniqueness of the solution in each case. This means that the design of mines is generally basic research in the case of a project or modernization of each mine. Adequate methods are required for the cooperation of specialists from various fields prepared for interdisciplinary communication. Achieving significant effects, if it is to be a result of rational and deliberate action, requires a solid team preparation or success in its construction.  This thesis is confirmed by the cooperation among  employees of Wrocław University of Technology on the design and analysis of open-cast mines that took place at the turn of the 1950s and 1960s. The leading motive of this cooperation was the reliability of the opencast mining engineering systems. However, it also had to run models of queueing theory and inventory control. A partial historical analysis of what was important in this matter can be found in Czaplicki's monograph \cite{Cza09:Shovel} (v. also \cite{Cza11:Meaning}).

Probabilistic methods in engineering sciences in Poland were applied by Witold Wierzbicki (1890-1965) (v. \cite{Jas1985:Wierzbicki}) in the lecture given on November 17, 1936, cf. \cite{Wie1936:ANT}). In the audience of this review there was a~mathematician, Witold Pogorzelski (v. \cite[p.9]{JozTer16:LODZ,Jer1988:Pogorzelski}), who participated in the discussion.

The presence  of these methods in the education of engineering majors dates to the early 1950s. For example, in the Department of Mathematical Statistics at Technical University of Łódź they appeared in 1954 (cf. \cite[p. 18]{JozTer16:LODZ}). Just after  World War II, there were only references to these methods in the syllabi of mathematical courses  (outline of probability calculus: Bernoulli theorem, Bayes formula; the theory of errors: Gauss's law, mean error and probable error.

It should be noted that in the curriculum of mathematics for engineering faculties of Wrocław University of Technology there could not have been more focus on probability and statistics. The courseload of mathematics was comparable to that of Łódź University of Technology, and the lecturers were mostly specialists in other areas of mathematics. Admittedly, prof. Steinhaus run a probabilistic seminar for students of mathematics, but that could not have significantly affected the education of engineers. Therefore, it should be assumed that the application of probabilistic models by Jan Sajkiewicz to mining engineering was due only to his own inquiries and the help of mathematician colleagues.


Kopocińska~\cite{KopIlona1968:Cyclic} applied the Takacs model (1962) to the analysis of the shovel-truck system operating in a gravel mine. The  model presented was of a two-stage type, and it was identical with M/G/I/m according to the Kendall notation. The problem of reliability was not addressed in the above papers. The problem of reliability, however, relates also to the organization of repairs and ensuring reserve systems. It is discussed in \cite[Ch. IV]{Kop1973:Outline}. In relation to mining engineering, special attention should be paid to \cite[\S\ 30, p.~272--282]{Kop1973:Outline} where the model of Sevastianov and Mayanovitch is presented.

Graff~\cite{Gra1971:Simple} discusses the problem of the application of the queuing model with an unlimited source of arrivals and times of an exponential character. 





\section{Mathematical models of open pit mines.}

In the years 1962-1999, first in the Chair of Mathematics, and later in the Institute of Mathematics and Theoretical Physics of  Wroclaw University of Technology\footnote{The Institute of Mathematics and Theoretical Physics, Wroclaw University of Technology was founded in 1968.}, over $20$ reports on the mathematical problems of exploitation of technological systems of open-cast mines with KTZ belt transport were made. Of fundamental significance were Stanisław Gładysz's works~\cite{Gla64:Process,Gla65:Wplyw,Gla65:Wydajnosc}. Technological and operational problems, as well as organizational problems related to the extraction of lignite and transferring it with conveyors to the power plant in the real functioning energy conglomerate in Turoszów exceeded the research and development capabilities of the Central Exploration and Design Mining Center (COBPGO)\footnote{Known as \emph{Poltegor}.}. The Union of Brown Coal Industries sent engineers and mathematicians of Wrocław universities packages of urgent problems to be solved. The vast number of tasks that should have been solved in the orders received was an inspiration to create a team composed of mathematicians and IT specialists solving scientific and application problems in close cooperation with COBPGO ,,Poltegor'', especially with prof. J.~Sajkiewicz (later the director of the Opencast Mining Institute of   Wrocław University of Technology). Over time, practical issues generated a~rich array of unresolved mathematical problems, especially those that were theoretically interesting, which became the main reason to expand the  team further by including younger employees, as well as professors from University of Wroclaw and their associates. In the project the following employees of Wroclaw University of Technology  were working: J.~Battek, R.~Kapala (Pomierski), E. Dorenfeld, E. Glibowski, E. Rychlikowski, T. Galanc, K.~Dyrka, K.~Janczewski, M.~Król and others (see reports \cite{DorGli66:WB,GlaBatSaj66:ProcesAwarii,BatGlaSaj73:Zarys,BatGlaSaj65:ProcesAwarii}). 


In the technological system Digger - Belt conveyor - Docker(DBD), two processes of decisive importance for the effective operation of the opencast mine were observed, namely the mass flow process (spoil in the mine, i.e. lignite or earth) $\{V(t); t \geq 0\}$ and the system of the failure process $ \{X (t); t \geq 0\} $.
\begin{figure}[th!]
\centerline{\includegraphics[width=0.76\textwidth]{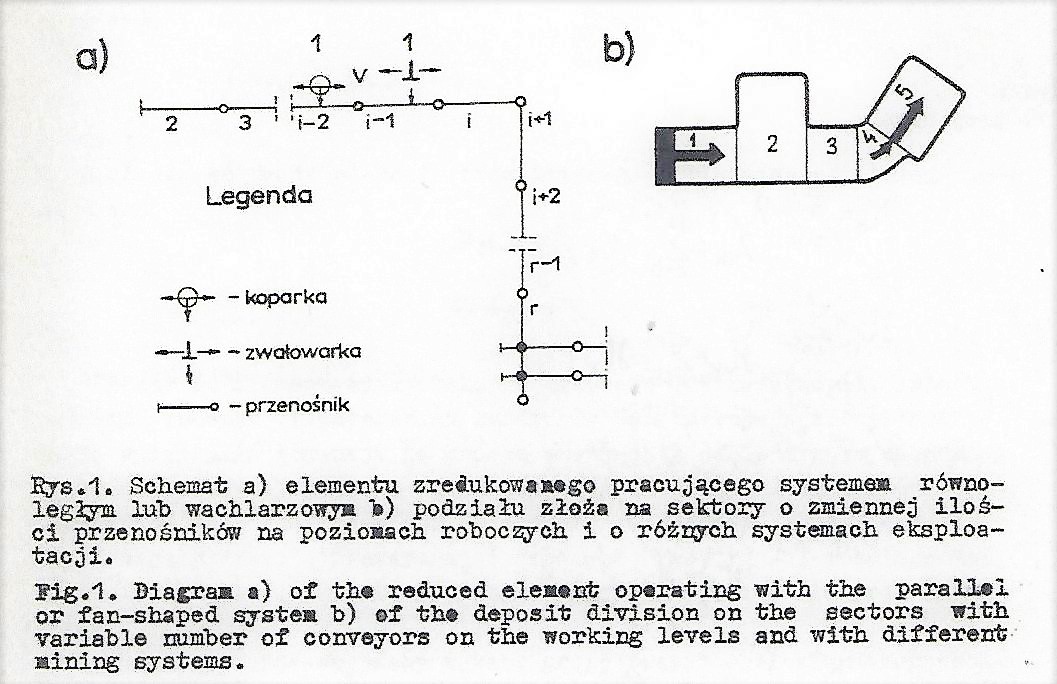}}

\centerline{\includegraphics[width=0.76\textwidth]{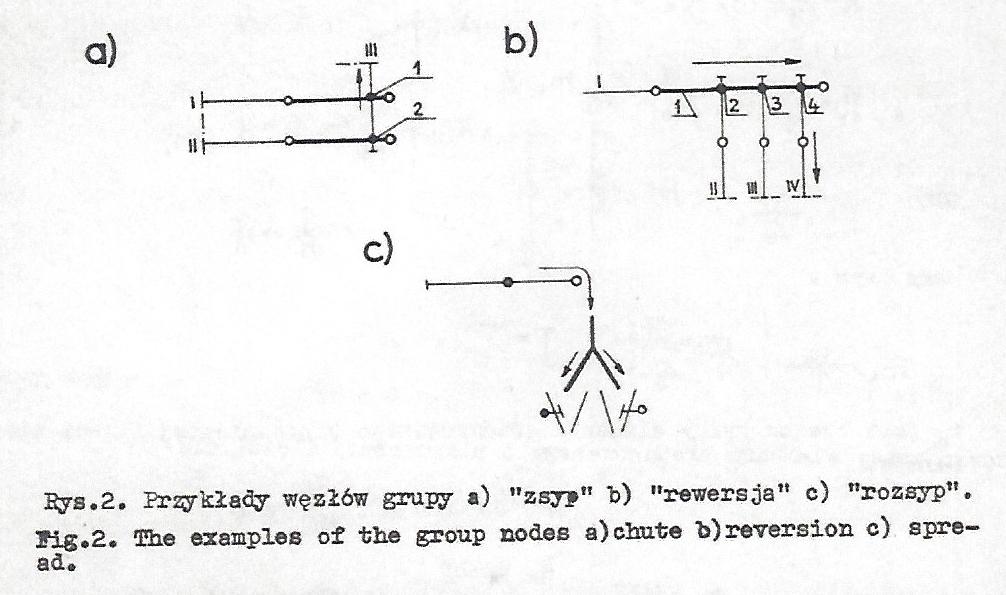}}

\footnotesize
\caption{\label{FigSchemeGladysz1}The key to the finalization of the model construction was the translation of detailed engineering schemes via  abstract graphs of the mathematical model  (v.\cite{BatGlaSaj73:Zarys}, \cite{Saj1979:Outline}). Important elements of the modelling were the group nodes.}
\end{figure}%
\normalsize

These processes described the phenomena determining the overflows on conveyor belts and overfill in reversing funnels. They also had a decisive impact on the size of reservoirs (bunkers), significantly improving the availability of DBD systems \cite{Gla65:Wydajnosc}. The basic components of DBD systems, namely excavators, conveyors and stackers, are complex devices in themselves. However, in the mathematical models created by the professor and his team, their meaning was elementary: at any given time they could only be in one of the following states: work, breakdown or forced stop. Usually the explained variable was the efficiency of the entire mining equipment system and for this reason stratification of the failure condition according to technological or organizational criteria did not have a logical justification \cite{Gla64:Process}. In a large proportion of the models created, excavator performance was determined by an appropriate weighted average only in mathematical considerations regarding the issues of calibrating the conveyor belts servicing an excavator or excavator systems; the efficiency was treated as a random variable whose probability distribution could change over time. According to the indications of prof. Gładysz,  rather than building  one global mathematical model describing the operation of the entire open-cast mine, the researchers focused on creating mathematical models of its simplified fragments, the simplification being strongly related to the purpose function and the issues to be investigated. In this way, the models were relatively transparent and allowed one to  answer  specific questions posed by engineers exploiting opencast mines \cite{Gla65:Wplyw}.

The greatest achievement of the team working on the problems arising from the mine's operation was the creation of mathematical models of the most important elements that make up the technological system.  For the implementation of such synergy, it is necessary to reliably translate technical problems contained in engineering projects in such a way that the abstract mathematical model is understandable for the constructors of the system. Due to such diligence in modeling processes and relationships between static parameters, engineers can describe the conclusions from an abstract model with real, measurable effects.  A good example of some communication problems are drawings in reports. Each important issue has an illustration in the form of a technological description with an illustration and presentation of the  related mathematical model with the scheme. This is an extremely important part of the project showing the determination of participants in adjusting the model in order to achieve the highest compliance with the real process modeled.

\section{\label{OCT} The effects of the involvement of mathematicians.}
A summary of the effects of the involvement of mathematicians in solving industrial problems in Lower Silesia is still waiting for conclusions. It is obvious that the  effort of a few people implies theoretical questions and the basic research in mathematics. The involvement of several people has shown significant opportunities to support project work with advanced mathematical modeling. The dynamics of industrial processes mean that the elementary methods used in designing are not sufficient to predict the behavior of complex technical systems. Attempts to describe such a system give the opportunity to evaluate the course of the exploitation and the preparation of alternative scenarios. All participants of the project made a significant contribution to its implementation. Knowledge of existing solutions and potential limitations of engineers forced mathematicians to look for a proper description. In turn, this forced the conduct of basic research, which resulted in, among other things, doctoral dissertations $x-y$. These dissertations arose from the study of topics and their solutions were expected by those who implemented the project. The small number of publications in mathematical journals related to the project is puzzling.

\begin{table}
	\centering\scriptsize 
	\setlength{\extrarowheight}{4pt}
		\begin{tabular}{|@{}c@{}|l@{}|p{17em}|l@{}|@{}c@{}|} \hline\hline
		\,No\,& The author&The title of disertation&The advisor&\,Year\,\\ \hline\hline
1&Adam Fiszer \cite{Fis1972:PhD}&The method of optimal planning of reservoir exploitation within the limits of the mining of the copper ore mines.&Jan Sajkiewicz&1972\\ \hline
2&Mieczysław Król~\cite{Kro1973:PhD}&About a certain variant of the linear stochastic process with barriers controlled by the Markov process.\protect\footnotemark&Stanisław Gładysz\,\, & 1973\\ \hline
3&Quang Nguyen \cite{Ngu1973:PhD}&Operational reliability of a continuous, double transportsystem with an intervention tank.&Stanisław Gładysz&1973\\ \hline		
4&Tadeusz Galanc \cite{Gal1973:PhD}&On one-dimensional distributions of a Markov, stationary, stochastic control process.&Stanisław Gładysz&1973\\ \hline
6&Radomir Simic \cite{Sim1974:PhD}&A method of selecting the optimal set of multi-bucket excavators for the exploitation of deposits.&Jan Sajkiewicz&1974\\ \hline
7&Jerzy Malewski \cite{Mal1974:PhD}&Transformation and synthesis of schemes of a certain class of machine systems with continuous transport and reliability of these systems.&Jan Sajkiewicz&1974\\ \hline
8&Hien Pham-Van \cite{Pha1974:PhD}&Influence of geometrical parameters of the open pit work on the selection of an economically justified depth of selecting thick hard coal deposits deposited in concise compositions for the case of loading discharged waste rocks with bucket excavators.&Jan Sajkiewicz&1974\\ \hline
9&Piotr Maćków \cite{Mac1975:PhD}&The method of testing the mining machinery operation process and assessment of their work efficiency.&Jan Sajkiewicz&1975\\ \hline
10&Witold Niedbala \cite{Nie1976:PhD}&Method for the identification of the operation process of machine systems with discontinuous technological structures.&Jan Sajkiewicz&1976\\ \hline
11&Zdzisław Iwulski \cite{Iwu1976:PhD}&The method of calculating technological firing parameters in carbonate rocks in LGOM mines.&Jan Sajkiewicz&1976\\ \hline
12&Krystyna Pazdyka \cite{Paz1976:PhD}\,&Analytical method for determining optimal mining and transport systems with discontinuous technological structures.&Jan Sajkiewicz&1976\\ \hline
13&Teresa Woźniak \cite{Woz1977:PhD}&Rheological model of machine operation for the control of material resources.&Jan Sajkiewicz&1977\\ \hline
14&Tomasz Dałkowski&Digital modeling of continuous transport systems.&Jerzy Battek&1977 \\ \hline
15&Ryszard Kabat \cite{Kab1978:PhD}&Digital modeling of subsystem operation consisting of a chain excavator on a rail chassis and objects cooperating with it at the working level.&Jan Sajkiewicz&1978\\ \hline
16& Stefan Sajkiewicz&Modeling of the exploitation of unloading mining and transport systems in mining.\protect\footnotemark&Jerzy Battek
&1980 \\	\hline 
\end{tabular}
	\caption{The dissertations related to the project.}
	\label{tab:TheDisertations}
\end{table}\normalsize
\footnotetext{The referees of the dissertation were: Bolesław Kopociński and Stanisław Trybuła.}
\footnotetext{The referees of the dissertation were: Michał Hebda, Jerzy Marcinkowski and Tadeusz Żur.}

Another important aspect of this project is its undoubted positive impact on strengthening the role of mathematics in the training of engineers. Initially, the concept assumed fuller education of engineering elites as part of an inter-faculty study (SPPT). Towards the end of the 1960s, the PWR Senate approved the transformation of the SPPT into the WPPT, which means that each PWR graduate was treated as a potential partner in projects that required advanced mathematical modeling techniques.

The importance of mathematical methods and exceptionally complicated procedures of their application have also motivated work on the availability and unification of the proven methods. Industry in Poland, following the example of Japan and other countries, strived  for quality control techniques. Mathematical methods are so complex that they give alternative outcomes. There is a need to compare them, indicate the right ones or decide if they are satisfactory. Systems of industrial standards describing the application of the best tested  methods of design, control and maintenance are  created. 

The history of the creation of the norms is a source of methodological inspiration and therefore requires special care and interest. 
\appendix
\section{Reports on ,,Życie Warszawy''.}

\subsection{\label{AppZW1964:Szperkowicz23vii1964}Mathematics on the outcrop. A billion we do not need to bother our heads about\protect\footnote{Jerzy Szperkowicz~\cite{ZW1964:Szperkowicz23vii1964}. Translated by Joanna Wasilewska.}.}  
Two poles of the topic -- the imagination of a mathematician from Wrocław and a mobile chute of a conveyor belt at the Turoszów mine. On the way there is a politician's office, designers' cardboard world and a checkered land of limits - the office of an investor.

I spent a few days in these realms.

I crossed the boundaries of concepts, which are clearer than borders at which a passport is requested.

I listened to explanations in languages of different continents of knowledge, understanding only the basic fact:

Stanisław Gładysz PhD, from the department of mathematics at the Wroclaw University of Technology, developed a method of programming conveyor belts for opencast mines. For sulfur in Machów, lignite in Turów and Konin and for new, future basins.

A bond between science and industry --  as my editor called it. The success of our center -- a Wroclaw politician will say. The result of developing basic research – a representative of a university, authorities will assess. Attracting science to life problems - but do not write too much – asks the head of the lignite projects office. Certainly, but $ \ldots $ - an investor will sigh tearfully. A solution to the given issue - the creator of the revealing method will say modestly.

I can hear an unexpectedly calm commentary in the place I was not advised to go at all – from engineers in the mine.

,,From the mathematical conclusions, we have chosen ten suggestions for the functioning and reconstruction of conveyor belts. We are introducing them gradually.''

Along the road from Turoszów to Bogatynia there is a range of gray hills covered with rare birchwood. Five years ago - I remember - the area on the opposite side of the road and housing estates. The road and houses remained. The hills wandered.

This is ultimately the basis of an open pit mine.

The excavators are biting into the ground's layers covering the deck, stackers - mechanical monsters - are forming hill ranges a few kilometers away. These are unknown to cartographers.

Immense masses of soil which are found here one can try to transport by trains, overhead cable-ways, but conveyors - rubber chutes with an electric drive - are increasingly being used. 

It is unbelievable how easily one could write about kilometers of conveyor belts, mountains to be moved, boilers of 20-storey power plants, only five years ago, all in future tense.

At that time I was assisting in covering with concrete foundations for the first cold store. I was guessing. Will it reach the passing jackdaw? Will the Miedzianka river accede to flow backwards? Who on the shuttle bus will find their life chance on the construction site and who will bury it here?

Then I found myself among material facts. With a power plant and cold stores obscuring a quarter of the horizon, with the rebellious Miedzianka river, that escaped underground, with Jerzy Bodziński M.Sc.Eng., who was initially a manual worker there, and with a deputy director of the investment department, who has recently lost the position of a clerk.

Forgive me, but I asked about mathematicians and conveyor belts only after the initial confrontation.

\begin{center}$\star$\end{center}
In a room of a hotel in Wrocław I try to translate an articulate lecture by the designer Sajkiewicz PhD into the language of a daily newspaper:\\
The appetite of a power plant is a known magnitude. In order to provide it with the necessary coal, you need to take off and transport definite masses of ground. One could use a sufficient number of excavators, stackers and build as many conveyor lines as needed to carry out the task in each situation:

But the economic calculation is staggering. One stacker more or one stacker less is the difference of 100 million zloty in  foreign currency. The difference of an imported factory.

The point is, therefore, to provide for the next about 40 years without an interruption the supply of coal for boilers, at the lowest possible cost, taking into account all possible failures and stoppages. 

Guided by the approximate calculations, the engineers designed the Tur\'ow system of excavators, conveyors and stackers.

Mathematicians calculated precisely the optimal system for  Turoszów's conditions.

Errors, which designers and practitioners did not avoid, were revealed in equations with mathematical accuracy. They will be gradually removed and avoided when designing next mines.

There are more conveyor belts working in Turoszów than in the rest of Europe. The next national specialty. For our own use and export. A~few years ago, we took our first lessons in the GDR.

\begin{center}$\star$\end{center}
In one of the cramped rooms I meet coworkers, Gładysz, Assoc. Prof., the giant Jerzy Battek\footnote{In original article his name is wrongly speeled ,,Batek''.}, PhD, and, younger than me, Glibowski and Galanc, MSc. Turoszów is distant history. They are now completing an urgent order for the Adamów mine.

- Was it difficult? - I risk the question.

- The most important - replies the Associate Professor - was the idea and finding the common determinant for such diverse calculation elements as the arrangement of reversing points (allowing to transfer cargo from one conveyor to another), the quality of propulsion engines, and even the level of labourers' technical culture.

- I understand - I answer with a purely rhetorical habit. I am also asking about the feeling of isolation, in the face of the inability to share, for example with those closest to me, what mathematics is passionate about, devouring, celebrating. A flash of agreement.

- It is about looking for some analogies in the material world. For example, dust particles vibrating in the light band give the image of certain concepts from the theory of probability.

- Recently, the interest of the press has arisen ...

- I agreed to talk to you at the request of our authorities. They reckoned it could be useful.

Our conversation took one and a half hour. We both considered it a~success.

\begin{center}$\star$\end{center}

- What our colleague Gładysz did, did not come instantly.

In 1951, professors Edward Marczewski and Hugo Steinhaus ran a~seminar on stochastic processes (the theory of probability). It was the direction of research leading to today's achievements. The issue of creating a renowned mathematical community - the Wroclaw school, which prides itself on the alumni following its beliefs, such as professor Kazimierz Urbanik, a high school graduate in 1948, the creator of the new sub-field in the probability calculus - generalized stochastic processes.

- And the working conditions- you have seen yourself. The Institute of Mathematics at the University does not have its own place. It's squatting at the Polytechnic.

\begin{center}$\star$\end{center}
Professor Bolesław Iwaszkiewicz, the ,,father'' of the city of Wroclaw, by education - as is known - a mathematician, recently returned from a~trip to the United States. Lots of impressions, but home worries are stronger:

- The industry clearly withdrew from the initial grand recognition of support for mathematicians.

In reply to this, Kumorkiewicz, an engineer from the Coal Lignite Union says:\\
- We have built the Turoszów Basin. The basin full of prototypes. Trial excavators, stackers and conveyor belts keep the largest Polish power plant in motion. Mathematicians found errors in our solutions. It is rather flattering for us, as it is a confrontation of reality between construction and the theoretical calculus. One talks about building certain ineffective devices and the possibility of saving - in future - a billion zlotys. Let us assume that we will cut spending on conveyor belts in the next mines. In different conditions, there may be other, more expensive issues than in Turów. It's easy to say - we'll save a billion.

\begin{center}$\star$\end{center}
It was not so much  working out a mathematical method, but rather its financial reflection that moved politicians, harnessed scientists, their wings and frightened investors. Mathematical programming of the systems - an excavator - a conveyor belt -  promises to reduce their costs by 10 percent. On the basis of plans we can calculate 16 billion for this purpose and 10\% of this  equals to one billion six hundred million. Approximately - a billion.

In the splendor of all this, hopes for improving the condition of Wroclaw's science were revived, they could be heard in the speeches of politicians and beams of light shone at the highest levels. Here, however, they remained cold blooded. A billion saved, bravo - applause - one can transfer some of the money from opencast mines and spend on something else.

This billion is indeed not suitable yet for the use in financial operations, but in the long run - and yet we do not live until tomorrow in this country - we do not need – I will use our Warsaw dialect - to bother our heads about it.

\subsection{\label{AppZW1964:Poza11xii1964}Mathematics outside the Outcrop (to the editor of ,,Życie'').}\footnote{Jerzy Szperkowicz~\cite{ZW1964:Poza11xii1964}. Translation: Joanna Wasilewska.}
After the feature ,,Mathematics on the Outcrop'' (,,Życie Warszawy'' July 23, this year) about the application -- by the Wroclaw mathematicians -- mathematical methods for the analysis of mining machines and conveyor belts in the Turów lignite mine, we received a letter from the Investment Bank in Warsaw that included, among others, the following statements:
\begin{enumerate}\itemsep2pt \parskip0pt \parsep0pt
\item[I] W. Bojarski, PhD, Eng. from the Investment Bank developed a~mathematical method for the analysis of mining machinery and conveyor belts, operations earlier than Stanisław Gładysz, PhD, Assoc. Prof., and his cooworkers from  Wroclaw University of Technology.
\item[II] Gładysz, PhD, Assoc. Prof., was a reviewer of the work by W.~Bojarski, PhD, Eng., dedicated to this issue.
\item[III] Having at his disposal a mathematical machine, Gładysz, PhD, Assoc. Prof., specified in a certain way and developed the previously described method (by Bojarski PhD, Eng.).
\end{enumerate}
Following the principles of good journalistic practice, we sent a copy of the letter from the Investment Bank to Gładysz, PhD, Assoc. Prof., asking to take a stand on the issue that was becoming delicate.

Coincidentally, Gładysz, PhD, Assoc. Prof., was not in the country at that time. Hence a delay in his  response, which in turn did not reach in the country the author of the feature ,,Mathematics on the Outcrop'', to whom the response by Gładysz, Ph.D, Assoc. Prof., was addressed.

Meanwhile, on November 26th this year ,,Życie'' released the letter from the Investment Bank, without printing the explanations by Gładysz PhD, Assoc. Prof., which was our aim when undertaking multilateral correspondence with the interested parties. Today we include extensive excerpts from the reply by Gładysz, PhD, Assoc. Prof., apologizing for depriving him of an opportunity to express his stand, for which we asked ourselves. We hope that this will at least partially rectify our negligence. Now, regarding the allegations of the Investment Bank:
,,I~believe – as Gładysz PhD, Assoc. Prof., says - that it is only a coincidental, though very significant misunderstanding.''
Digital machines, arithmometers or slide rules are not used and have never been used to refine or develop methods. In my study, they do not play any role to the extent that in any of my publications there is  not even a mention of digital machines. Digital machines and similar devices are  used only to accelerate numerical calculations, which occasionally might be very significant. But the method must be prepared in advance ...
Each author of scientific publications needs to accept responsibility to familiarize themselves with the current state of the situation, the results of other employees and other centers. Bojarski PhD, was neither the first nor the last author who used probability calculus methods for these and similar issues. Simply put, the probabilistic language is here the natural language ...

The method which was developed by myself and continued by J.~Batek and J. Sajkiewicz is not only a ,,refinement'' and development of the Investment Bank method, but it is something completely different in the qualitative sense. It operates on entirely different concepts and, most of all, treats phenomena dynamically - in time - as they are really occurring. The omission of the basic time factor, as it takes place in all the studies known  to me, including that by the Investment Bank, makes it impossible to properly describe the work of the technological system. The point is that time cannot be somehow attached to this or that pattern, but one has to create a new theory from the scratch. This is one of the main difficulties successfully solved at this moment.

It is difficult to get engaged in more details here. However, I would like to draw attention to the fact that my results had already been presented many times to specialists, and on many occasions this happened  in public. Many of them were perfectly familiar with the development of the Investment Bank, but none of them had ever found any analogies.

And finally, a modest proposal. Soon the IVth National Mining Congress ,,New Mining'' will be held in Katowice. I was invited by organizers together with J. Battek and J. Sajkiewicz to give the lecture ,,The stochastic model of the working system -  excavator - conveyor belt - stacker''. I suggest that the Investment Bank comes forward and gives a presentation as well. The direct confrontation of both studies in front of the probably most competent group of specialists should allow to explain the misunderstanding. 

With the last post we received a letter from the Rector of the Wroclaw University of Technology, professor Z. Szparkowski. The letter expresses surprise and regret that the replica of the Investment Bank had been printed without the inclusion of Gładysz PhD Assoc. Prof., explanations, for which the reason is given above. As to the substance, professor Szparkowski states in the letter:
,,The Institute of Mathematics has appointed a special committee, which after hearing both the work of Bojarski PhD., Eng., as well as the work of Gładysz, PhD, Assoc. Prof., stated that they are completely incomparable ... ''

A comprehensive explanation following the issue was sent by the Lower Silesian Mining Project Bureau located in Wroclaw.

At the last minute, we also received a letter and a statement from representatives of the Wroclaw mathematical community. They conclude that the letter by the Investment Bank published in ,,Życie'' did great harm to Stanisław Gładysz, PhD, Assoc. Prof., who is the author of the currently used (and the only appropriate) method of calculating the effective efficiency of complex transport systems in opencast mines – and they demand rectification of the letter of November 26th.

The letter and the statement were signed by: Edward Marczewski - full professor, member correspondent of the Polish Academy of Sciences, head of the Institute of Mathematics at University of Wroclaw, Adam Rybarski - Associate Professor, head of the Mathematics Department at  Wroclaw University of Technology, Hugo Steinhaus - full professor, full member of the Polish Academy of Sciences, Władysław Ślebodziński - full professor, former head of the Mathematics Department at  Wroclaw University of Technology, Kazimierz Urbanik - full professor, head of the Probability Calculus Department and vice-rector at University of Wroclaw.

Once again, we express our regret referring to the oversight of simultaneous publications of stands on the feature ,,Mathematics on the Outcrop''. As for the substance of the controversy, for understandable reasons, we believe that if there is anything else to add, a more appropriate place would be the forum  mentioned at the end of the letter by Gładysz, PhD, Assoc. Prof.

\section{Biographical research}
\subsection{\label{JerzyBattek}Jerzy Battek  (1927--1991)}
\begin{wrapfigure}{l}{3.3cm}
  \vspace{-8mm}
  \begin{center}
    \includegraphics[width=3.1cm]{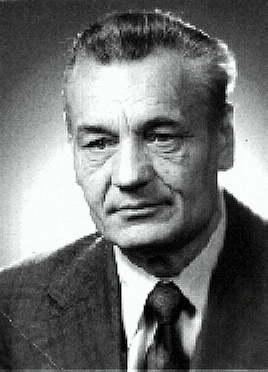}
  \end{center}
  \vspace{-20pt}
\end{wrapfigure}
\begin{figure}[th!]
\centerline{\includegraphics[width=90mm]{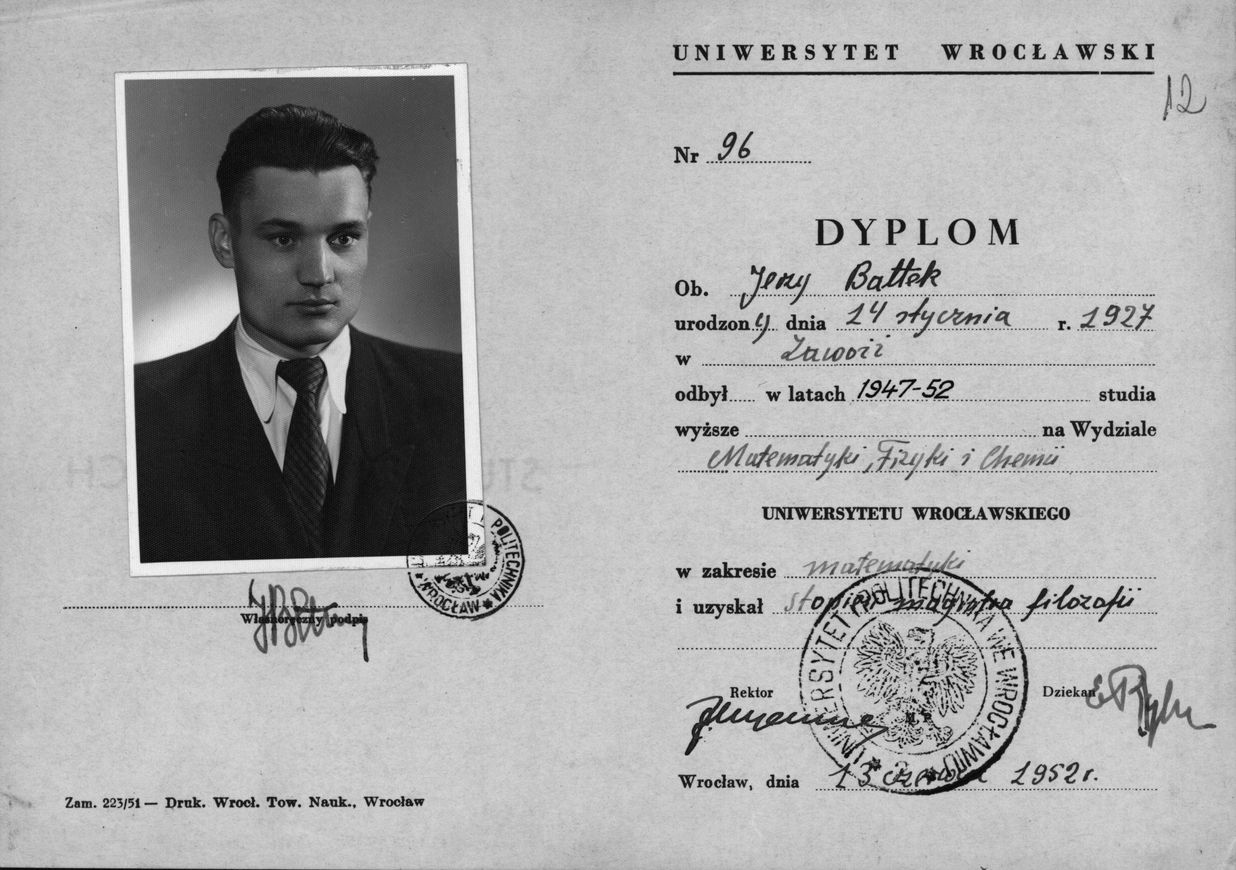}}
\caption{\label{JBDyp1}Diploma of completing higher studies with a master's degree in philosophy (\cite[mks 22 AUWr]{EGlibUWr1963}).}
\end{figure}

was born on January 14, 1927 in Zawoja\footnote{The bio of Jerzy Battek (1927--1991) can be found in  \cite{Huz2006:Battek} (v. \cite[p.89]{DudWer06:WMS}).} . He started his education in elementary school and graduated from his last class in Krak\'ow. The outbreak of war prevented him from attending middle school. He spent the war in Krak\'ow, attending clandestine lessons, thanks to which in 1945 he passed  exams for the small high school diploma. After the war, he went to Lower Silesia and Continued his education at the State High School for Adults in Wroclaw. After the matriculation exam in 1947, he entered the Faculty of Mathematics, Physics and Chemistry at University of Wrocław, from which he graduated in 1952 with a master's degree in philosophy (v. Fig.~\ref{JBDyp1}).

He started his professional work while studying, working in the years 1949--1952 as a teacher in secondary education in Wołów and Wroclaw. After graduating, he started working as a deputy assistant at the Department of Mathematics at Wroclaw University of Technology, headed by Professor Władysław Ślebodziński at the time.
 Jerzy Battek was the founder and the first director of the Computation Center, the predecessor of the IT Center, with the title of Docent dr.  His tall, slim and slightly stooped figure and inseparable pipe were known to everyone who came into contact with the Computation Center. And there were a lot of these people, because the period in which he came to work was a pioneering period in the use of computers at Wroclaw University of Technology. His name is inseparably connected, in the opinion of the Wroclaw University of Technology, with the creation and functioning of the Computing Center as the first university-wide unit offering computer services. He was among  pioneers who laid technical and organizational foundations for future applications of computer science at Wroclaw University of Technology. Battek was the supervisor of two doctoral dissertations (see Table\,\ref{tab:TheDisertations})\,directly\,related\,to\,the various\,topics\,in\,maintenance\,of\,mines.

\subsection{\label{AKErlang}Agner Krarup Erlang(1878--1929)} was born on 1st of  January 1878 in Lønborg, near Tarm, in Jutland. He was a Danish mathematician, statistician and engineer. He was a son of a schoolmaster, who was also a parish clerk. His mother was related to  Thomas Fincke (1561--1656), a Danish mathematician and physicist. In 1892, he passed the Preliminary Examination of the University of Copenhagen with distinction, after receiving dispensation to take it because he was younger than the usual minimum age. For the next two years was teaching alongside his father. He started studies in 1896 and graduated in 1901. After that he was teaching mathematics for seven years. He met the Chief Engineer of CTC\footnote{Copenhagen Telephone Company} at the meeting of the Danish Mathematical Society and from 1908 he began to work for this company. His close relation to CTC lasted twenty years until his death in Copenhagen after an abdominal operation on the 3rd of February 1929. He published his ideas in  Danish mathematical journals. His most important papers in queueing theory related to the traffic engineering were \cite{Erl1909:telefons,Erl1920:Telefon} (v. \ref{Erlang1920}).

\subsection{\label{EdmundGlibowski} Edmund Glibowski (1928--2014)} 
\begin{figure}[th!]
\centerline{\includegraphics[width=9cm]{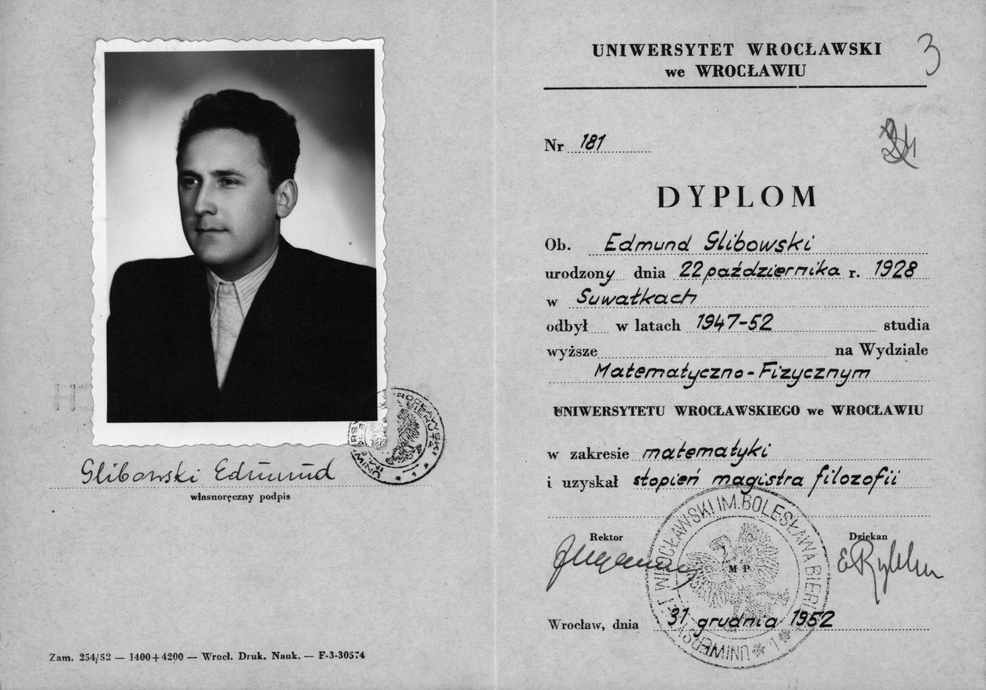}}
\caption{\label{EBDyp1}Diploma of completing higher studies with a master's degree in philosophy
 (\cite[mks 3 AUWr]{EGlibUWr1963}).}
\end{figure}

He was born on November 22, 1928 in Grodno\footnote{Currently a city in Lithuania.}{${}^{,\,}$}\footnote{Glibowski gave this birthplace in a CV submitted along with the application for the commencement of the doctoral proceedings\cite[mps 2]{EGlibUWr1963}. In other documents we can read for the place of birth: Suwałki, woj. białostockie.}, a son of Stefan and Julia ne\'e Boronowska. His father Stefan was a railwaymen (he was shot on March 13, 1943). He started his education in elementary school in Grodno. After the World War II, he came to Poland and stopped in Świętokrzyskie, where he attended a general high school. He graduated from high school in 1947 in Wroclaw, in High School No. 1. In this school he later worked as a mathematics teacher in the years 1950-1952. After graduation, he studied at the Mat-Fiz-Chem Faculty of University of Wroclaw. He graduated in 1952 (v. \ref{EBDyp1}) and took up a job at the Chair of Mathematics of Wroclaw University of Technology on September 1, 1952. 

\begin{figure}[th!]
\centerline{\includegraphics[width=9cm]{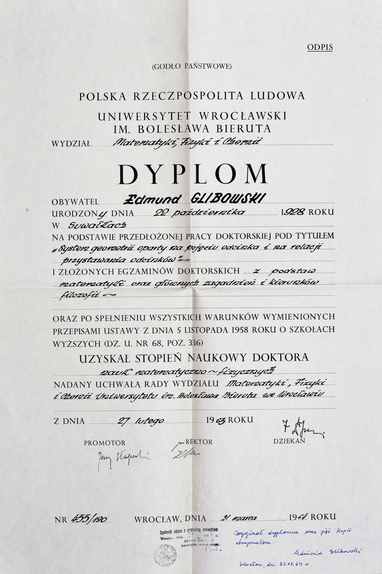}}
\caption{\label{EBDyp2}Edmund Glibowski's PhD Dilpoma.}
\end{figure}

On June 19, 1963, he applied for  opening of  doctoral proceedings with thesis entitled ,,Systems of elementary geometry based on the notion of a segment of the real line  and the relationship of congruent segments''\footnote{,,System geometrii elementarnej oparty na pojęciu odcinka i relacji przystawania odcinków''.}. Professor Jerzy Słupecki was appointed the supervisor of the doctoral dissertation. He passed his doctoral exams in the main fields of philosophy and the foundations of mathematics on February 27, 1963. On the same day in the afternoon he also defended his thesis\footnote{The announcement of defense appeared on February 8, 1963 in the ,,Słowo Polskie''. The defense was carried out at 5p.m. in the Anthropology Room at the Polish Academy of Sciences at ul. Kuźnicza 35 on the first floor in room 13. The course of the defense was recorded by prof. Stefan Paszkowski.}. The reviews were provided by Bronisław Knaster (July 25, 1962) and Lech Dubikajtis (November 21, 1962). The examining committee was headed by doc. dr Zygmunt Galasiewicz (physicist) and its members were Lech Dubikajtis and Jerzy Słupecki. The facsimile of diploma is presented in~\ref{EBDyp2}.

Starting from March 1, he was promoted to assistant professor. In addition to basic work in geometry, for which he obtained his doctoral degree, he keenly cooperated with engineers. At the time when ,,Życie Warszawy'' was interested in the activity of mathematicians in cooperation with engineers on the issues of the efficiency of an energy system based on fossil raw materials, he noted at least one report on this subject \cite{DorGli66:WB}. In addition, he was interested in the theory of graphs, which he used in the modeling of ceramic materials (in works with Zbigniew Święcicki).

In private life, he was married twice. He retired on October 31, 1988 due to poor health. He died on March 16, 2014.

\subsection{\label{StGladysz} Stanisław Gładysz (1920–2001)} was born on March 22, 1920 in Piotrów near Opatów, Świętokrzyskie Voivodeship, son of Józef and Emilia ne\'e Zimmerman. He attended a general high school in Kraków and graduated in mathematics in 1948 from the Jagiellonian University (v. \cite[mks. 7]{StGlaUWr1961}). He worked at the Gliwice Polytechnic (1956--48), and from 1948 at Wroclaw Polytechnic till 1985 and then in 1986--1991. In 1956 he was appointed for a position at the Mathematical Institute of the Polish Academy of Sciences. 
\begin{figure}[th!]
\ifwww
\centerline{\includegraphics[width=10cm]{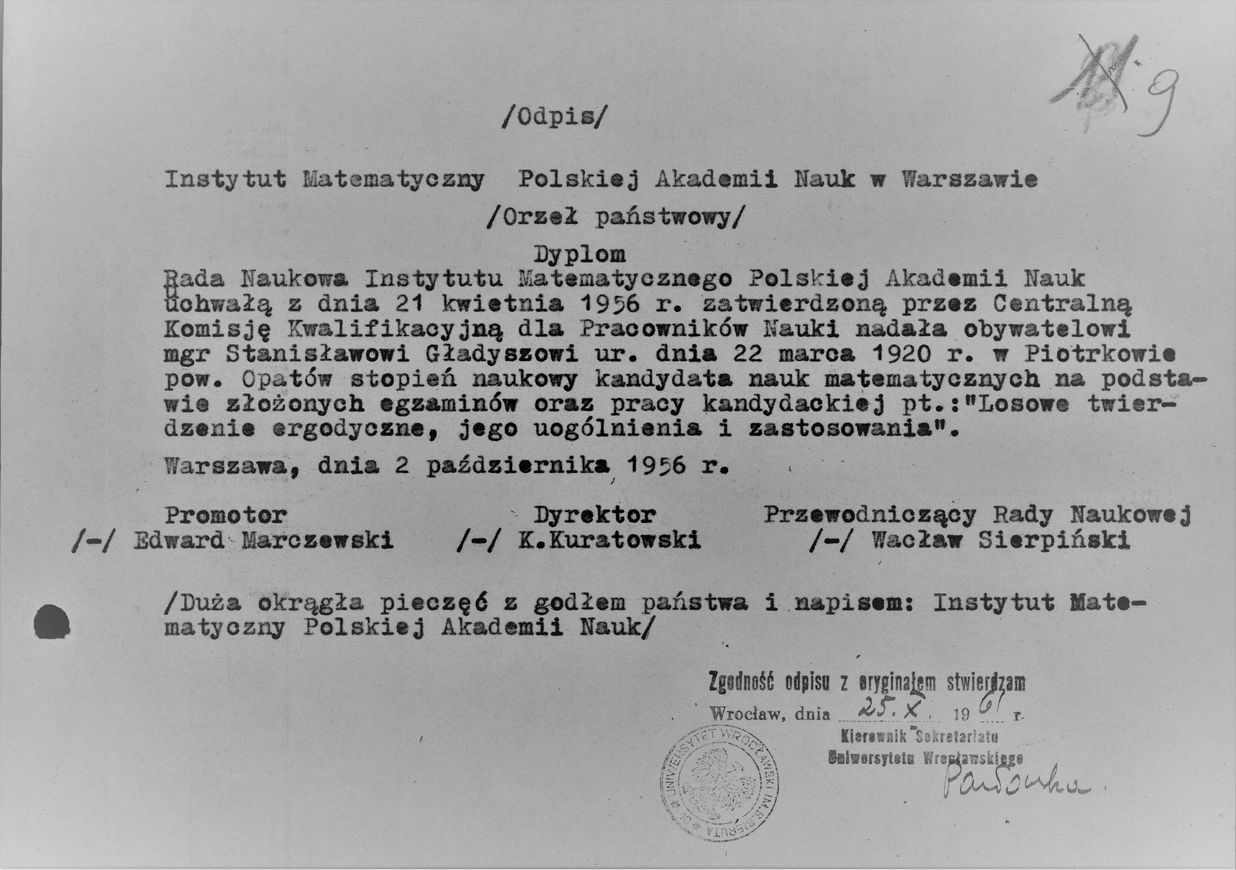}}
\else
\centerline{\includegraphics[width=10cm]{GladyszStPhD.jpg}}
\fi
\caption{\label{GlaDoc1}The diploma with the description of the doctoral thesis (mks AUWr syg. IV-4010/1961 Gładysz Stanisław poz. 9).}
\end{figure}

\begin{figure}[th!]
\centerline{\includegraphics[width=105mm]{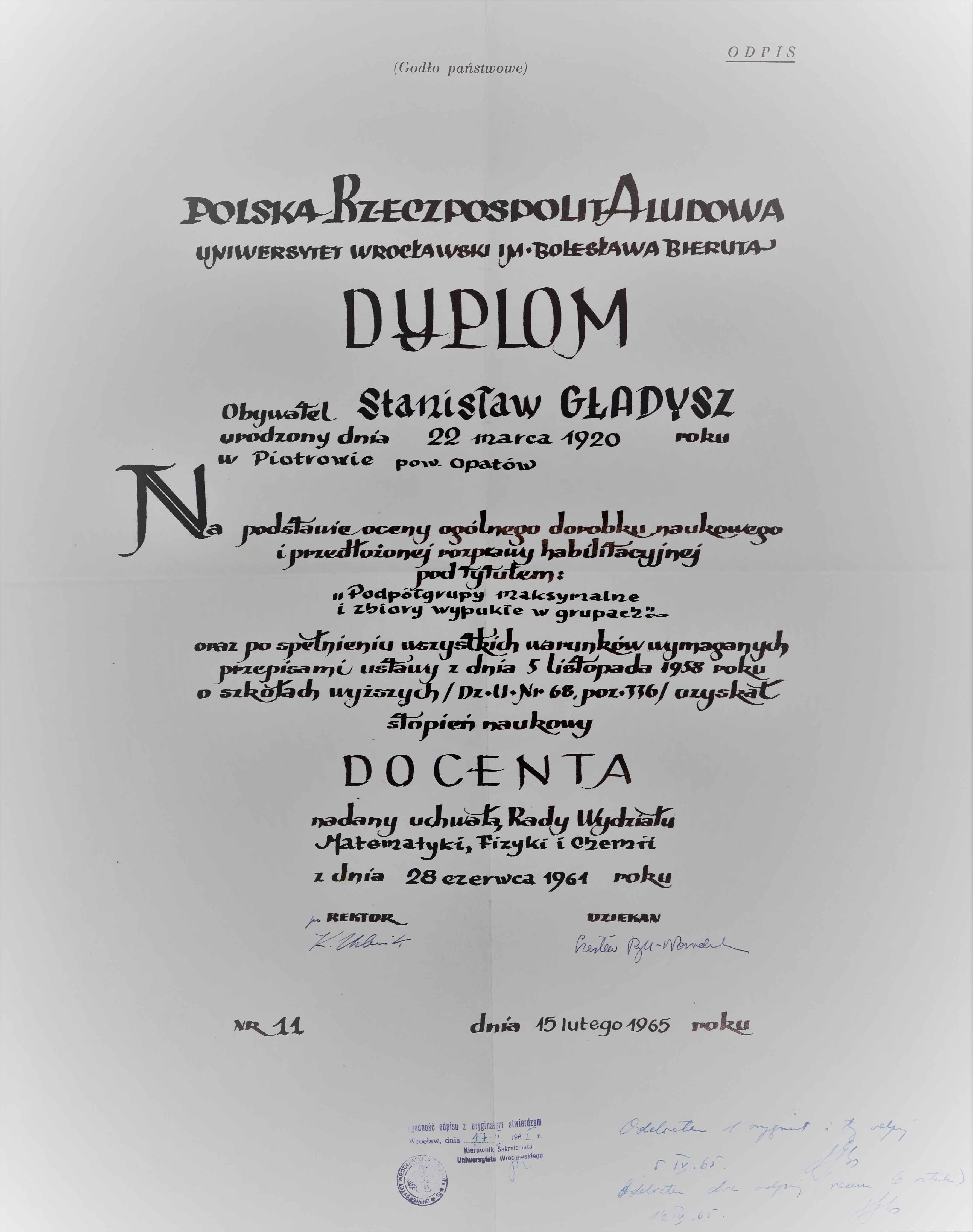}}
\caption{\label{GlaDoc2}The Assistant Professor Diploma (mks AUWr syg. IV-4010/1961 Gładysz Stanisław poz. 40).}
\end{figure}

On April 21, 1956 he defended the doctoral dissertation ,,A random ergodic theorem, its generalization and applications'' at the Institute of Mathematics of the Polish Academy of Sciences in Warsaw (v. Fig.~\ref{GlaDoc1}). Professor Edward Marczewski was appointed the supervisor of the dissertation and the opponents were Stanisław Hartman PhD, Prof., Czesław Ryll-Nardzewski PhD, Prof. and, as the C.K.K reviewer, Edward Marczewski PhD, Prof. On April 11, 1961 he applied for opening of the procedure of habilitation at Mat-Fiz-Chem. of Wrocław University based on the dissertation ,,Maximum semigroups and convex sets in groups'' (v. \ref{GlaDoc2}). The reviewers were professors Czesław Ryll-Nardzewski and Kazimierz Urbanik. The opinions about the candidate were delivered by professors: Roman Sikorski, Helena Rasiowa and Stefan Zubrzycki. He worked at Wroclaw University of Technology  in the years 1949--1985, and  1986--1991, occupying the following positions: deputy professor (since 1955), assistant professor (since 1965) and professor (since 1967).
 
This was his period of creative scientific work and implementation of the bold concept of creating an independent Institute of Mathematics at Wroclaw University of Technology. In addition to valuable theoretical results related to the subject of his thesis, he had significant achievements in cooperation with scientific units dealing with other disciplines, as well as with industry. As far as theoretical works of Stanisław Gładysz are concerned, among the most important are works on random versions of ergodic theorems and interesting reformulations of ergodic theorems in the language of ergodic functionals, introduced axiomatically. The professor also studied properties of the space of measurable sets with the metric introduced by Hugo Steinhaus, as well as the properties of the maximum subgroups that were used to  define abstractly a convex topology in groups. Because of his research interests   he lectured willingly for many years on topology and elements of harmonic analysis on LCA groups at Wroclaw University of Technology and his experience as an academic teacher is included in the  original lecture notes and a~topology textbook. A separate category of works concerned stochastic processes and their applications (cf. the section~\ref{RelModHis}). Their  value among other results is due to the comprehensive solution to the reliability problem of (excavator-conveyor-dumping system--OCT) BDB system obtained by Stanisław Gładysz and a team of mathematicians in cooperation with the Mining Institute and Poltegor. It was a very successful attempt at  application of Markov processes and ergodic theory, which brought measurable economic benefits to the energy conglomerate in Turoszów.
 
The turning point in the practical implementation of Gładysz's plans and visions was the period of 1969--1974. These years are characterized by significant structural changes in the employment (first graduates of the School of Fundamental Problems of Technology study at the Institute of Mathematics and Theoretical Physics, the Faculties of Civil Engineering, Electronics and Opencast Mining). The Institute received the right to award the doctoral degree and  the first doctoral dissertations were defended. During this period, prof. Gładysz was the supervisor of five doctoral dissertations directly or indirectly related to this research topic. The defense of four of them (T. Byczkowski, T. Galanc, M. Król, A. Weron) took place at the Institute of Mathematics and Theoretical Physics, and only one, of a~trainee from Vietnam Nguyen Quanga, at the Mining Institute of Wroclaw University of Technology.

\subsection{\label{JanSajkiewicz}Jan Sajkiewicz (1928--1995).} 
\begin{figure}[th!]
\ifwww
\centerline{\centerline{\includegraphics[width=10cm]{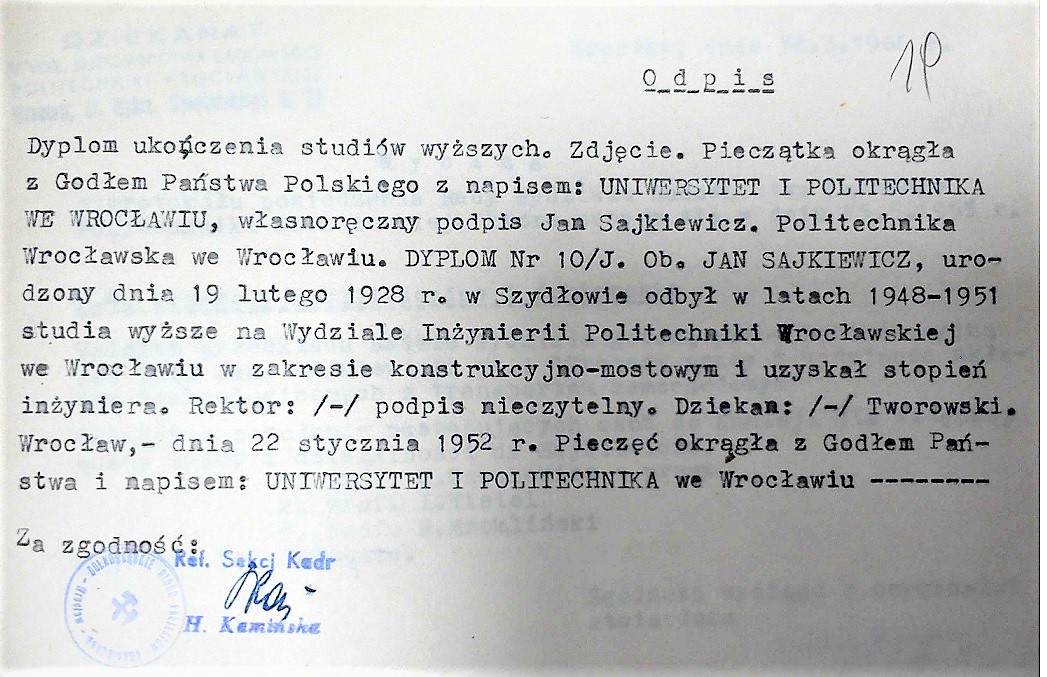}}
}
\else
\centerline{\centerline{\includegraphics[width=10cm]{OdpisDypInzJSaj1952a.png}}
}
\fi
\caption{\label{SajDoc1}A copy of the diploma of higher education with the degree of an engineer (mks APWr id 163/1 15-D-207 folio 19).}
\end{figure}

Jan Sajkiewicz, the son of Feliks and Bronisława ne\'e Trelińska, was born on 19 February 1928 in Szydłów (Staszów county) in a family of farmers (v. \cite{JanSajAPWr1928}).  He lived in the native village until 1932, when the whole family emigrated to Argentina for profits. They returned a year later. He graduated from the elementary school in Szydłów in 1940. In winter 1942, he worked out with his (two years older) brother the secondary school program, and passed an exam on two years of high school program  in April 1944. After the Red Army entered  Szydłów in 1944 and withdrew fast (after 2 weeks),  he stayed with his parents on the farm. His brother joined the Second Polish Army in the meanwhile. In 1947, he graduated from Jan Śniadecki's  high school (in Kielce) in Mathematics and Physics program. After graduation, he worked for a year as the headmaster of the Primary School in Gacki (Szydłów commune). In 1948, he began studies at the Faculty of Engineering of the Wrocław University of Technology and on 22 January 1952 he received an engineering diploma in the field of construction and bridges (A copy of the Diploma \ref{SajDoc1} from \cite{JanSajAPWrHab1928}).


\begin{figure}[th!]
\centerline{\includegraphics[width=105mm]{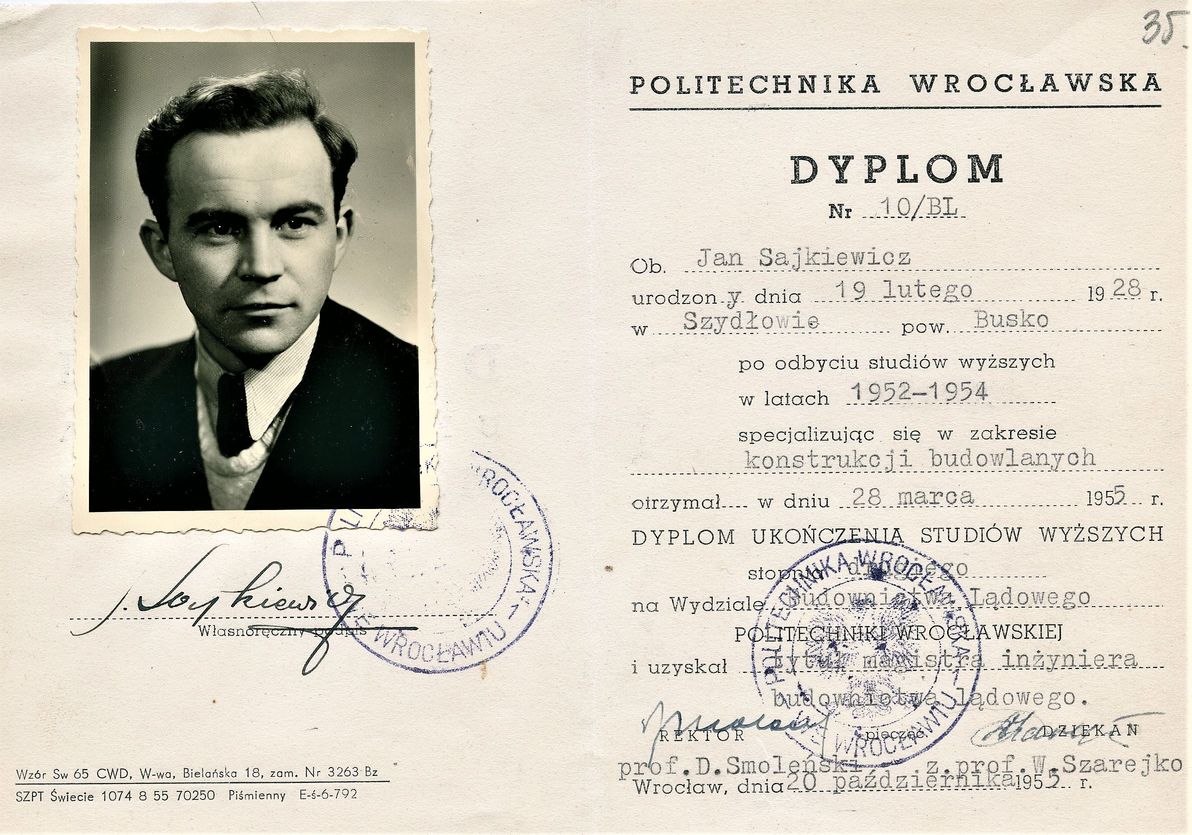}}
\caption{\label{SajDoc2}Diploma of university graduation with the degree of engineer (mks APWr id 163/1 15-D-207 folio 35).}
\end{figure}

In the years 1950--1952 he held a scholarship  and was junior assistant of the Chair of Construction Structure at the Faculty of Engineering at Wroclaw University of Technology. Then he worked in the army as a researcher and at the same time studied at the Faculty of Civil Engineering at Wroclaw University of Technology (v. copy of the obtained diploma \ref{SajDoc2}), from which he  graduated in 1955. At that time he specialized in the field of building constructions. From 1956, he worked in the Lower Silesia Mining Project Office(\emph{Poltegor}) as the head of the studio, then the general designer of the ,,Gosławice'' mine, ,,Adamów'' mine and the mines of the Turek, Lodz and Konin region. He also lectured on the economics of the organization and planning of opencast mining.


He got the building construction license in 1959 (according to Article 362 of the Construction Law). He obtained the degree of doctor of technical sciences defending the dissertation ,,\emph{Determination of the maximum technically possible extraction of mineral useful from the opencast mine}'' on May 25, 1963. The supervisor of the dissertation was Prof. Bolesław Krupiński (v. mks APWr, id 217/86 folio 8).  

\begin{figure}[ht!]
\centerline{\includegraphics[width=90mm]{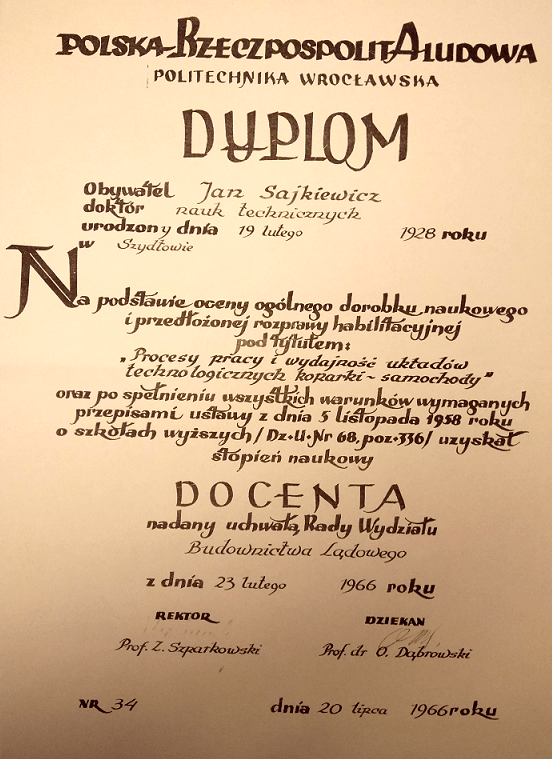}}
\caption{\label{SajDyplDoc}Diploma awarding the scientific degree of Associate Professor to Jan Sajkiewicz PhD.}
\end{figure}

He defended the habilitation in 1966 (v. \cite{JanSajAPWrHab1928}), and obtained professor's title on July 6, 1973 (v. \cite{JanSajAPWrTytul1928}). 

He was employed as the  deputy director of the research department at ZBiP Miedzi ,,Cuprum'' , 1968--1970, and then as the director. He was the director of Institute of Mining of Wroclaw University of Technology in 1971--1981 (after Prof. W. Czechowicz) and at the same time the director of the Operational and Prognostic Research Department. In 1978, he became honorary doctor of Mining Academy in Freiberg. In the years 1981--1982 he was the Vice-Rector for Didactics and Education at Wroclaw University of Technology.

He became the Rector of the Higher School of Engineering (WSI) in Radom from 1982--1985 and the chairman of the Department of Machine Operation and Technical Facilities at Kielce University of Technology in 1985, and he chaired in 1991--1993 the Department of Agricultural Machines and Equipment of the same Technical University.

He is the author of over 300 items: books (issued in Poland and abroad), lecture notes, articles, and has a share in solutions protected by patents. He supervised 18 doctoral dissertations in engineering, as a specialist in the field of operation and reliability of facilities and technical systems recognized in Poland and abroad. His numerous achievements were prized and won awards.
\selectlanguage{english}

\subsection{\label{EHWWeibull}Ernst Hjalmar Waloddi Weibull (1887--1979)} was born on 18 June 1887.  He was a Swedish engineer, scientist, and mathematician. His biography intertwines with the challenges of the twentieth century. He joined the Swedish Coast Guard in 1904 as a midshipman. He moved up the ranks with promotion to sublieutenant in 1907, Captain in 1916 and Major in 1940. As the coast guard officer he had an opportunity to take courses at the Royal Institute of Technology. In 1924 he graduated and became a~full professor \cite{Bro97:Weibull}. Weibull died on October 12, 1979 in Annecy, France (cf. \cite{Bro97:Weibull}, \cite{Stochastikon2017:EHWWeibull}).

\section{Markov models in standarization issue.} Technical systems designed by engineers should, in addition to fulfilling the expected tasks, also ensure the highest quality of work, increase the reliability of equipment and employees safety, etc. Therefore, there has always been a need for some routine preventive activities to guarantee the expected effect. For this purpose, the International Organization for Standardization (ISO) was established in London in 1947, which proposes standards to Member States. The Polish Committee for Standardization (PKN) was one of the founding countries and is a permanent member of this organization. PKN, which was established in 1924, issued its first standard in 1925.

Standards are  based not only on practice, but also on mathematical theory and its correct applications. Quantitative methods based on probability theory are a recognized tool for their creation. However,  establishing  rules for more technologically advanced processes  requires the use of models based on Markov processes. Their description, using a graph of transitions between reliability states, is so clear and effective that it has been an approved basic tool for technical sciences. This method of modeling has been incorporated into standard methods and described in the standards to provide comparable assessments of the issues modeled. The first standard concerning the application of Markov's processes in Poland was PN-IEC 1165: 1998 (see \cite{NormaPN-EN61165})) (current version of PN-EN 61165: 2006). The methods contained there are based on the same principles as those in Gładysz's works in, e.g. \cite{Gla65:Wplyw}.

\medskip

\noindent{\bf Acknowledgements.} The creation of this study was prompted by Dr. Rościsław Rabczuk, who mentioned the role of journalists in popularizing applied mathematics in the times of the Polish People's Republic. Undoubtedly, discussions with dr hab. Marian Hotlos on the occasion of establishing the beginnings of the Faculty of Fundamental Probems of Technology was inspiring. The friendly help of Mrs. Anna Rubin-Sieradzka from the \href{http://www.archiwum.pwr.edu.pl/index.php}{University Archive of the Wrocław University of Science and Technology} was invaluable. Without the encouragement of Professor Domoradzki and a warm word from Professor Monika Hardygóra, the Dean of the Faculty of Geoengineering, Mining and Geology, Wrocław University of Science and Technology, it would not be possible to reach the stage of submitting the article for publication. The authors would like to thank all the mentioned and some other people who contributed to the creation of this study. 
\bigskip

\noindent{\large\bf References}\nopagebreak

\medskip
\selectlanguage{polish}

\begin{center}
{\bf  Projektowanie systemów inżynierskich w polskich kopalniach w drugiej połowie XX wieku}\\[1.5ex]
\href{\repo/5126
}{Aneta Antkowiak, Monika Kaczmarz \& Krzysztof Szajowski}
\end{center}

\begin{abstract} Udział matematyków w realizacji projektów gospodarczych w Polsce, w których realizacji metody matematyczne odegrały istotną rolę, pojawiał się sporadycznie w przeszłości. Metody, zwykle sprawdzone i znane z opublikowanych opracowań, są adaptowane do rozwiązywania pokrewnych problemów. Przedmiotem prezentacji będzie współpraca jaką nawiązali matematycy i inżynierowie we Wrocławiu w drugiej połowie XX wieku przy analizie efektywności systemów inżynieryjnych stosowanych w górnictwie. Efekty tej współpracy spowodowały wprowadzenie metod projektowania systemów technicznych wcześniej nieznanych inżynierom. Otworzyły też tematy badawcze w~teorii procesów stochastycznych i teorii grafów. Interesujące są też społeczne aspekty tej współpracy.     
\end{abstract}

\selectlanguage{english}

\noindent Aneta Antkowiak\\ 
\textit{E-mail address}: \textbf{221721@student.pwr.edu.pl}\\[3ex]
Monika Kaczmarz\\ 
\textit{E-mail address}: \textbf{Monika.Kaczmarz@pwr.edu.pl}\\[3ex]
Krzysztof Szajowski\\ 
Wrocław University of Science and Technology\\
Faculty of Pure and Applied Mathematics, \\
Wybrzeze Wyspianskiego 27, PL-50-370 Wrocław, Poland\\
\textit{E-mail address}: \textbf{Krzysztof.Szajowski@pwr.edu.pl}

\end{document}